\documentclass[preprint,11pt,authoryear]{elsarticle}
\usepackage[margin=4.3cm]{geometry}

\usepackage{amssymb}
\usepackage{amsmath}
\usepackage{amsmath,amssymb,amsfonts}
\usepackage{algorithm, algpseudocode}
\usepackage[mathscr]{euscript}
\usepackage{multicol}
\usepackage{multirow}
\usepackage{xcolor}

\begin{document}

\begin{frontmatter}



\title{Distributionally Fair Peer-to-Peer Electricity Trading}

\author[inst1]{Estibalitz Ruiz Irusta}
\ead{estibalitz.ruiz@uni.lu}
\affiliation[inst1]{organization={University of Luxembourg},
            addressline={29 Av. John F. Kennedy}, 
            city={Luxembourg},
            postcode={1359}, 
            country={Luxembourg}}

\author[label1,label2]{Juan Miguel Morales}
\ead{juan.morales@uma.es}
\affiliation[label1]{organization={Group OASYS},
            addressline={Ada Byron Research Building, Arquitecto Francisco Peñalosa Street, 18}, 
            city={Málaga},
            postcode={29010}, 
            country={Spain}}
\affiliation[label2]{organization={Statistics and Operations Research, Universidad de Málaga},
            city={Málaga},
            country={Spain}}
\begin{abstract}
Peer-to-peer energy trading platforms enable direct electricity exchanges between peers who belong to the same energy community. In a semi-decentralized system, a community manager adheres to grid restrictions while optimizing social welfare. However, with no further supervision, some peers can be discriminated against from participating in the electricity trades. To solve this issue, this paper proposes an optimization-based mechanism to enable distributionally fair peer-to-peer electricity trading. For the implementation of our mechanism, peers are grouped by energy poverty level. The proposed model aims to redistribute the electricity trades to minimize the maximum Wasserstein distance among the transaction distributions linked to the groups while limiting the sacrifice level with a predefined parameter. We demonstrate the effectiveness of our proposal using the IEEE 33-bus distribution grid, simulating an energy community with 1600 peers. Results indicate that up to 70.1\% of unfairness can be eliminated by using our proposed model, even achieving a full elimination when including a non-profit community photovoltaic plant. 
\end{abstract}


\end{frontmatter}

\section{Introduction} \label{introduction}

Power systems continue evolving from a hierarchical organization towards a decentralized and user-centric structure \citep{1_energy_collectiveness}. While investment costs decrease, the deployment of Distributed Energy Resources (DER), such as solar and wind, is growing worldwide \citep{2_fair_and_effective}. As a result, end-use customers are now becoming prosumers, as they no longer consume energy only, but also produce \citep{3_emerging_buss}. Consequently, end-use customers are transitioning from being passive consumers to having active roles within the power system operation. 

To facilitate this transformation and ensure a secure operation of the power system, new markets and mechanisms are needed \citep{4_social_vulnerability}. Peer-to-peer (P2P) electricity markets allow the direct exchange of electricity without the need for intermediaries \citep{5_resilience}. Peers can buy electricity from other peers with a surplus. Such a decentralized market brings peers the possibility to obtain economic benefits while increasing the transparency and efficiency of the power system \citep{6_block_enabled}. Furthermore, it can empower participants to trade electricity locally and promote renewable energy integration \citep{7_leveraged}. 

Different settings exist within a P2P market structure. When having a fully decentralized market configuration, peers trade electricity directly without the supervision of a third party \citep{6_block_enabled}. However, the presence of a supervisory role allows the consideration of network constraints to avoid congestion issues, as well as the optimization of resources to increase social welfare. Therefore, by incorporating a community manager to the problem, peers can communicate with a central node that supervises the exchange to help in the decision-making process \citep{1_energy_collectiveness}. 

Nevertheless, when we grant an algorithm the authority to decide for citizens, concerns about fairness arise \citep{8_review}. Researchers already started analyzing the disparities between peers in P2P electricity markets, mentioning, for instance, that peers who cannot afford generation technologies will never sell electricity \citep{sharehold}. Additionally, peers who live in isolated or marginal areas will more likely have lower distribution capacities, limiting their possibilities to trade electricity. Consequently, it is essential to consider social justice aspects when designing a P2P electricity market.

The notion of \emph{fairness} and \emph{equity} hold different meanings across the literature \citep{4_social_vulnerability}. Some authors focus on determining a fair trading price within the P2P market. We can identify two main approaches in this direction. The first perspective advocates for a uniform trading price for all participants, determined by the intersection of supply and demand curves \citep{2_fair_and_effective}, \citep{7_leveraged}. However, since peers often operate under different tariff agreements with utilities, this uniform pricing may not reflect individual circumstances. To address this, authors in (\citep{5_resilience}) propose a model in which each peer submits a buying or selling bid specifying their minimum or maximum acceptable price. Then, an optimization algorithm creates as many matches as possible to maximize the social welfare of the community, while defining a different trading price per match.

Another central view of justice is the fair allocation of goods \citep{8_review}, where most authors see this good as the earned revenue \citep{fair_profit, renewable}. Back in ancient Greece, Aristotle already introduced the idea of distributive justice. He argued that goods should be allocated to individuals based on what they rightfully deserve \citep{8_review}. Since then, numerous philosophers and intellectuals have attempted to define what each person deserves. 

Nonetheless, previous research has not analyzed what exactly causes unfairness within a P2P electricity market. Moreover, it focuses entirely on \emph{individual} and not \emph{collective} fairness. Thus, even when fairness actions are introduced, socially vulnerable groups can still be discriminated against as a community. To solve these issues and fill the associated research gaps, this paper focuses on analyzing the sources of unfairness within P2P electricity markets. Additionally, we propose a distributionally fair optimization model to minimize the disparities among different groups defined based on their energy poverty level. More precisely, we define a measurable index to quantify unfairness based on the \emph{amount of electricity traded} and use this index to analyze which variables can impact fairness within P2P electricity markets. We then propose a distributionally fair optimization-based mechanism to enhance equity in P2P electricity markets using the Wasserstein metric. Our mechanism requires solving an NP-hard problem, which we address using an alternating algorithm. Finally, we demonstrate the effectiveness of our proposal using the IEEE 33-bus distribution grid with 1600 peers.

In Section \ref{method} of this paper, we describe the details of our methodological approach. Afterward, in Section \ref{solution}, we clarify the adopted solution strategy. Then, Section \ref{case_study} presents the case study employed for validation, while the results are discussed in Section \ref{results}. Finally, Section \ref{conclusions} concludes the paper by summarizing our work and offering future recommendations. 

\section{A distributionally fair P2P electricity trading mechanism} \label{method}

Consider a community represented by a set $\mathcal{P}$ of peers willing to participate in a P2P electricity market. Each peer $i \in  \mathcal{P}$ has a consumption $c_i$ and electricity production $e_i$ at a certain point in time. When $c_i \leq e_i$, this peer has a surplus to sell, and thus, sends a bid specifying the minimum price at which she is willing to sell each kWh of surplus ($\lambda^{s}_i$). Conversely, when $c_i \geq e_i$, this peer has a deficit and sends a buying bid specifying the maximum price at which she is ready to buy each kWh of deficit ($\lambda^{b}_i$).

The maximum and minimum prices at which each peer is willing to buy and sell electricity are directly linked to the prices contracted and agreed with the utility. For instance, if the utility is selling electricity at $\lambda^{us}_i$ to peer $i$, then this peer buys electricity from another peer only if the offering price is lower than or equal to $\lambda^{us}_i$. Similarly, if the utility is willing to buy electricity from peer $i$ at a price $\lambda^{ub}_i$, then the minimum price at which this peer is willing to sell her surplus will be higher than or equal to $\lambda^{ub}_i$. 

In a fully ``selfish'' P2P market setting, the aim is to maximize the sellers' revenue, as expressed in \eqref{eq:ref_obj}, mimicking the so-called ``sell to the highest bidder'' system. In  \eqref{eq:ref_obj}, \(X\) is the transaction matrix deciding which peer sells to who and buys from whom. Each element $X_{ij}$ reveals how much electricity (kWh) peer $i$ sells to peer $j$. Notice that the real revenue of participating in the P2P market is not the direct profit of the transactions, but the extra profit made in comparison with trading with the utility. Additionally, we establish that the final price of a transaction will be the average of the selling and buying bids of the peers who match. Nonetheless, not all transactions are possible. There exist certain constraints that the community manager must respect.

\begin{subequations} \label{eq:REF}
    \begin{align}
        &\underset{U_i^{buy}, U_i^{sell}, X_, e_i^{c}}{\max}\quad \sum_{i \in \mathcal{P}} \sum_{j \in \mathcal{P}} \left( X_{ij} \left(\dfrac{\lambda^s_i + \lambda^b_j}{2} - \lambda^{ub}_i\right)  \right) \label{eq:ref_obj} \\
        s.t.\quad &  e_i - e_i^{c} - c_i  = \sum_{j \in \mathcal{P}} X_{ij} + U^{buy}_i - \sum_{j \in \mathcal{P}} X_{ji} - U^{sell}_i  &\forall i \in \mathcal{P} \label{eq:ref_balance1}\\
        & e_i^{c} \leq e_i &\forall i \in \mathcal{P} \label{eq:ref_e_i_c}\\
        & X_{ij} \leq Y_{ij} M_i \quad & \forall i,j \in \mathcal{P} \label{eq:ref_lim} \\
        & \underline{v} \leq 1 + 2A^{-1}D_rA^{-T} p_n + 2A^{-1}D_xA^{-T}q_n \quad & \forall n \in \mathcal{N} \label{eq:ref_bus_d} \\
        &  1 + 2A^{-1}D_rA^{-T} p_n + 2A^{-1}D_xA^{-T}q_n \leq \overline{v} & \forall n \in \mathcal{N} \label{eq:ref_bus_u} \\
        & p_{n} = \sum_{i \in \mathcal{P}(n)}{e_i- e_i^{c} - c_{i}} \quad & \forall n \in \mathcal{N} \label{eq:ref_p}\\
        & q_{n} = \sum_{i \in \mathcal{P}(n)}{\rho_{i}} \quad & \forall n \in \mathcal{N} \label{eq:ref_q}\\
        & \label{eq:ref_non-neg} U_i^{buy}, U_i^{sell}, X_{ij}, e_i^{c} \geq 0 \quad  & \forall i,j \in \mathcal{P}
    \end{align}
\end{subequations}

Firstly, the energy balance must be guaranteed within the P2P market \eqref{eq:ref_balance1}. Peers cannot buy more electricity than the total surplus of the community from other peers. If more energy is needed, the extra electricity should come from the utility.  Variables $U_i^{buy}$ and $U_i^{sell}$ represent the electricity bought from and sold to the utility by peer $i$, respectively. Moreover, variable $e_i^{c}$ offers curtailment options in case the distribution grid cannot handle the original production values. Additionally, if peers $i$ and $j$, $i\neq j$, do not match their bids because $\lambda^{s}_i >  \lambda^{b}_j$, it means that it is economically more beneficial to trade with the utility. Therefore,  \eqref{eq:ref_lim} ensures that only legitimate trades happen. By using  $M_i$, $i \in \mathcal{P}$, a big enough value, and \(Y\), a bid matching matrix containing 0 where a transaction is infeasible and 1 when it is feasible, we ensure bid matching. 

Selecting appropriate values for $M_i$, $\forall i \in \mathcal{P}$, is important to ensure model tractability and solution quality. If the constants $M_i$ are excessively large, it can result in numerical issues, which may hinder solver performance. Consequently, each $M_i$ must be ideally calibrated to be the smallest value that remains sufficiently large to avoid the exclusion of any feasible solution. Within our problem, the transaction matrix values are inherently bounded by the maximum potential surplus available to each peer, as no peer can trade more electricity than they generate. This generation capacity is constrained by the installed capacity of each peer's PV plant. Therefore, we define $M_i$ as the maximum installed PV capacity of peer $i$ in the network.

Finally, the community manager must guarantee that the distribution line capacities are not surpassed. We use the Linear Distribution Flow (LinDistFlow) model approximation to ensure such limitations \citep{LinDistFlow}. In a radial distribution network with $\mathcal{N}=\{0,\ldots, n\}$ buses and $\mathcal{L}=\{1,\ldots,l\}$ lines, where $|\mathcal{L}| = |\mathcal{N}|-1$, each bus $n \in \mathcal{N}$ has a voltage $V_n$. Additionally, the net power injection at bus $n$ is $S_n = p_n + jq_n$, where $p_n$ is the active power injection and $q_n$ is the reactive one. The LinDistFlow approximation linearizes the distribution flow equations by neglecting the line loss terms \citep{lindistflow2}. Equations \eqref{eq:LDF1}--\eqref{eq:LDF3} linearly relate the squared voltage magnitudes ($v_n=|V_n|^2$), the power injections ($p_n, q_n$), and the line flows ($P_n, Q_n$), where \(r_{(\pi_n, n)}\) and \(x_{(\pi_n, n)}\) are the resistance and reactance of the line connecting bus $n$ with its parent $\pi_n$, respectively. Moreover, \(v_n\) and \(v_{\pi_n}\) are the squared voltages of bus $n$ and its parent \(\pi_n\).
\begin{subequations}\label{eq:LDF}
    \begin{align}
        &\sum_{k:n \rightarrow k} P_k \thickapprox p_n + P_n, & \forall n \in N \label{eq:LDF1}\\
        &\sum_{k:n \rightarrow k} Q_k \thickapprox q_n + Q_n, & \forall n \in N \label{eq:LDF2}\\
        &v_n \thickapprox v_{\pi_n} - 2 (r_{(\pi_n, n)} P_n + x_{(\pi_n, n)} Q_n),\quad\quad\quad\quad\quad\quad &  \forall n \in N \label{eq:LDF3}
    \end{align}
\end{subequations}

By incorporating the reduced branch-bus incidence matrix $A$ with dimensions $|\mathcal{L}|\times|\mathcal{L}|$ we can write the matrix form of the LinDistFlow approximation, see Equations \eqref{eq:LDF_matrix}. Matrix $A$ describes the connections between the buses (excluding the substation node) and branches of the system and $a^\top_0 = (1, 0, \ldots, 0)$  is the column vector of length $|\mathcal{L}|$ associated with the substation bus. On the other hand, $D_r$ is a diagonal matrix with the elements of $r= (r_{(\pi_n, n)}, \forall n \in \mathcal{N}\setminus{0})$ as its diagonal entries \eqref{eq:diag}, while all off-diagonal entries are zero. Similarly, $D_x$ is the diagonal matrix corresponding to the reactance $x = (x_{(\pi_n, n)}, \forall n \in \mathcal{N}\setminus{0})$. Moreover, $p=[p_1,\ldots, p_n]^T$ and $q=[q_1,\ldots, q_n]^T$ are the active and reactive power injection vectors, and $P=[P_1,\ldots, P_n]^T$ and $Q=[Q_1,\ldots, Q_n]^T$ are the active and reactive power flow vectors, respectively. Additionally, $v=[|V_1|^2,\ldots,|V_n|^2]^T$ represent the square voltage magnitude vector, and $v_0$ is the squared voltage of the substation bus. 
\begin{subequations}\label{eq:LDF_matrix}
    \begin{align}
        & D_r = diag(r),  \quad D_x = diag(x), \label{eq:diag}\\
        & p = A^T P,  \quad q = A^T Q,\\
        & Av + v_0 a_0 = 2D_rP + 2DxQ,\\
        & v = v_0 + 2A^{-1}D_rA^{-T}p + 2A^{-1}D_xA^{-T}q, \quad\quad\quad\quad\quad\quad\quad\quad\quad
    \end{align}
\end{subequations}

Equations \eqref{eq:ref_bus_d} and \eqref{eq:ref_bus_u} ensure that the squared voltages of each bus $n \in \mathcal{N}$ are restricted within their secured ranges in per-unit-value ($\overline{v}, \underline{v}$) using the matrix form of the LinDistFlow approximation, where $\mathcal{P}(n)$ is the subset of the set $\mathcal{P}$ of peers that are connected to node $n$, such that $\mathcal{P} = \bigcup_{n\in N}\mathcal{P}(n)$. In this way, \eqref{eq:ref_p} and \eqref{eq:ref_q} calculate the active and reactive powers injected at or withdrawn from node $n$ in per-unit value. Peers with no PV installation will primarily consume reactive power $\rho_i$ due to the inductive nature of their loads. In contrast, the inverters of PV systems can both produce and consume reactive power, depending on their configuration and operating conditions. Therefore, $\rho_i$ can take both positive and negative values for each peer $i$ with an inverter. It is worth mentioning that the financial trade and the physical flows of electricity in a power system are detached. This means that a peer can still trade with another peer even if there is little or no physical connection between them, as long as they match their bids and the total power flows respect the limitations imposed by the distribution network.

\subsection{Identification of revenue}

Each household can make revenue by trading electricity with their peers and obtaining better conditions than those offered by the utility. By knowing the amount of traded electricity (kWh) stored in transaction matrix \(X\), the settled price (€/kWh), and the price offered by the utility (\texteuro/kWh), we can calculate the (extra) profit of each peer (\(\gamma_i\)):
\begin{equation}
\label{eq:revenue} \gamma_i  =  \sum_{j \in \mathcal{P}}\left( X_{ij}  \left(\dfrac{\lambda^s_i + \lambda^b_j}{2} - \lambda^{ub}_i\right) + X_{ji}  \left(\lambda^{us}_j -\dfrac{\lambda^s_j + \lambda^b_i}{2}\right) \right), \enskip \forall i \in \mathcal{P}\\
\end{equation}

\subsection{Identification of energy community groups}

Considering that our problem concerns power systems and the traded good is electricity, we believe that energy poverty is the appropriate metric to identify different groups. \emph{Energy poverty} is defined as the inability to secure materially and socially necessitated energy services, such as heating the house or using appliances \citep{poverty_assessment}. It has a multidimensional nature and is an indicator of domestic energy deprivation and energy vulnerability \citep{poverty_towards}. The energy poverty level of a household can be identified with the following variables: income level, the energy efficiency level of the house, available infrastructure, energy prices, and their volatility, health, age, and gender \citep{poverty_aristondo,poverty_revisiting,poverty_literature}. Using these variables, we can identify the energy poverty level of each peer and classify them into different groups. Peers can thus be organized into $k+1$ distinct groups, represented by the set $\overline{G} = \{g_1, g_2, \ldots, g_k, g_{pv}\}$. Each group $g_i, i \in \{1, 2, \ldots, k, pv\}$, is a contiguous subset of peers, such that $g_1=\{1, \ldots, m_1\}$, $g_2=\{m_1 +1, \ldots, m_2\}$, until $g_k=\{m_{k-1}+1, \ldots, m_{k}\}$ and $g_{pv}=\{m_{k}+1, \ldots, |\mathcal{P}|\}$. Importantly, the groups are mutually exclusive and collectively exhaustive, meaning every peer belongs to exactly one group $g \cap g' = \emptyset, \forall g, g' \in \overline{G}: g \neq g'$. This partitioning ensures a clear and non-overlapping organization of peers.

We remark that the set $\overline{G}$ includes a special group, $g_{pv}$, to represent actors excluded from the quantification of fairness. These actors participate in electricity trades but differ fundamentally from standard peers. In our model, we specifically consider non-profit community PV plants as such actors. These PV plants sell their electricity generation $e$ to peers but do not buy electricity, as they have no consumption $c$. Unlike standard peers, whose objective is profit-making, non-profit community PV plants aim to foster fairness within the energy community. To this end, they offer their electricity generation at no cost (i.e. $\lambda^s_i = 0, \forall i \in \{m_{k}+1, \ldots, |\mathcal{P}|\}$). 

\subsection{Unfairness level}

First, we need to determine the energy poverty level of each household. Once this information is obtained, we can assign each household to one of the groups $g$ from the set $G = \overline{G}\setminus\{g_{pv}\}$. Then, we need to calculate the total traded electricity, either buy or sell (kWh), that each peer $i$ made within each group $g \in G$. Using the information stored in the transaction matrix $X$ we can build one transaction distribution \(T^g\) for each group \(g \in G\) by summing all sold and bought electricity from each peer and actor $j \in \mathcal{P}$ \eqref{eq:transaction}.
\begin{equation}
    \label{eq:transaction}T_i^g =  \sum_{j \in \mathcal{P}} \quad (X_{ij} + X_{ji}), \quad \forall i \in g, \quad \forall g \in G \\
\end{equation}
The ``difference" between the distributions determines the unfairness level. We use the Wasserstein metric to calculate the distance between these distributions \eqref{eq:was1}--\eqref{eq:was2}. We can write this distance as an optimization problem that seeks to find the minimum transportation cost from one distribution to another, where $\pi^{g, g'}$ is the transportation plan, and \(d_{ij}^{g, g'}\) is the absolute distance from each support point $T^{g}_{i}$ to each support point $T^{g'}_{j}$ of the transaction distributions corresponding to groups $g$ and $g'$, respectively \eqref{eq:distance}. 
\begin{equation}
    d^{g, g'}_{ij} = | T^{g}_{i} - T^{g'}_{j} |, \quad \forall i \in g, j \in g', \quad \forall g, g' \in G: g \neq g'\\ \label{eq:distance}
\end{equation}
\begin{subequations} 
    \begin{align}
        \mathcal{W}(g,g')=\underset{\pi_{ij}}{\min}  & \sum_{i \in g}\sum_{j \in g'} \pi_{ij}^{g, g'}  d_{ij}^{g, g'} \label{eq:was1}\\ 
        s.t.\quad &  \sum_{j \in g'} \pi_{ij}^{g, g'} = \frac{1}{|g|}, & \forall i \in g \\
        & \label{eq:was2}\sum_{i \in g} \pi_{ij}^{g, g'} = \frac{1}{|g'|}, & \forall j \in g'
    \end{align}
\end{subequations}
The unfairness level (\(D_{max}\)) of the P2P model will then be the maximum distance from all the Wasserstein distances \eqref{eq:dmax}. 
\begin{equation}
    D_{max} = max\{\mathcal{W}(g,g'), \quad  (g, g') \in \mathcal{C}(G,2)\} \label{eq:dmax}\\
\end{equation}
where $\mathcal{C}(G,2)$ represents the 2-combinations of the set $G$. 
\subsection{Distributionally Fair Peer-to-Peer model}

The original optimization problem seeks a transaction matrix $X$ that maximizes the total profit of the sellers. However, it does not consider any aspect of fairness. To correct this, we now bring a sensitivity attribute into our problem, which is specified as the energy poverty level of the peers. This attribute will define the sensitivity set $G = \{g_1,\dots, g_k\}$ and partition the outcome space into different groups. Now, the goal of our optimization model will no longer be the maximization of revenue, but the minimization of the unfairness level \eqref{eq:fair_obj}. More concretely, we will aim to minimize the maximum Wasserstein distance between the transaction distributions of the various groups while maintaining the community profit around a near-optimal region \citep{Hanasusanto}.  

\begin{subequations} 
\footnotesize
    \begin{align}
        & \underset{ \mathcal{V}}\min \hspace{34mm} D_{max} && \label{eq:fair_obj}\\ 
        &s.t. \enskip (\ref{eq:ref_balance1})-(\ref{eq:ref_non-neg}) && \label{eq:fair_repeated}\\ 
        &  \sum_{i \in g} \sum_{j \in \mathcal{P}}   X_{ij} \left(\dfrac{\lambda^s_i + \lambda^b_j}{2} - \lambda^{ub}_i\right) +  X_{ji} &\hspace{-7mm}\left( \lambda^{us}_j -  \dfrac{\lambda^s_i + \lambda^b_j}{2} \right)  \geq (1-\epsilon)\cdot |\gamma_g^*|,\quad \forall g \in G& \label{fair_margin}\\
        & \sum_{i \in \mathcal{P}} U^{sell*}_{i} \geq  \sum_{i \in \mathcal{P}} U_{i}^{sell}  &&\label{fair_utility}\\
        & \sum_{i \in \mathcal{P}} e^{c*}_{i} \geq  \sum_{i \in \mathcal{P}} e_{i}^{c}  &&\label{curtail_e}\\
        & \sum_{j \in \mathcal{P}}  (X_{ij} + X_{ji}) = T^g_{i},  &\quad \forall i \in g, \quad \forall g \in G &\label{fair_transactions1}\\
        & t^{g,g'}_{ij} = T^{g}_{i} - T^{g'}_{j}, \quad s^{g,g'}_{ij} = T^{g'}_{i} - T^{g}_{j}, &\quad \forall i\in g,\quad \forall j \in g', \quad \forall (g,g') \in \mathcal{C}(G,2)& \label{eq:d1}\\
        & t^{g,g'}_{ij} \leq d_{ij}^{g,g'}, \quad s^{g,g'}_{ij} \leq d_{ij}^{g,g'},  &\quad \forall i\in g,\quad \forall j \in g', \quad \forall (g,g') \in \mathcal{C}(G,2)& \label{eq:d2}\\
        & \sum_{i \in g} \sum_{j \in g'} ( \pi_{ij}^{g,g'} d_{ij}^{g,g'} ) \leq D_{max},   &\quad  \forall (g,g') \in \mathcal{C}(G,2)& \label{eq:wd11}\\
        & \sum_{j \in g'} \pi_{ij}^{g,g'} = \frac{1}{|g|},  &\quad \forall i \in g, \quad \forall (g,g') \in \mathcal{C}(G,2)& \\
        & \label{eq:wd33}\sum_{i \in g} \pi^{g,g'}_{ij} = \frac{1}{|g'|},   &\quad \forall j \in g', \quad \forall (g,g') \in \mathcal{C}(G,2)& 
    \end{align}
\end{subequations} 

The optimization model is subject to the set of variables $\mathcal{V}$ = \{$D_{max}$, $U_i^{buy}$, $U_i^{sell}$, $X_{ij}, \forall i,j \in \mathcal{P}$; $T^{g}_{i}, \forall i \in g, \forall g \in G$; $t^{g,g'}_{ij}$, $s^{g,g'}_{ij}$, $d_{ij}^{g,g'}$, $\pi_{ij}^{g,g'}$, $\forall i \in g, \forall j \in g', \forall (g,g') \in \mathcal{C}(G,2)\}$. Equations \eqref{eq:fair_repeated} are the constraints from the reference optimization model where we ensure the energy balance, bid matching, line constraints, and the variables' non-negativity character. Additionally, equation \eqref{fair_margin} ensures that the collective profits made by each group do not deviate more than (\(1- \epsilon\)) from the original profits, where \(|\gamma_g^*|\) is the original profit made by each group $g$ \(\in\) $G$. Equation \eqref{curtail_e} ensures curtailment does not exceed the reference level ($e^{c*}_{i}$) while allowing to curtail from a different peer $i$ to enhance fairness. On the other hand, equation \eqref{fair_utility} ensures that fairness is not achieved by trading more with the utility and taking from the community, but only by distributing the trades better within the community. Afterward, to enable the Wasserstein distance calculation, we must compute the distance from each transaction $i$ from distribution $g$ to another transaction $j$ in distribution $g'$ for all possible pairs of groups $(g,g')$ chosen from $G$. To calculate the absolute distance value, we first introduce auxiliary variables $t_{ij}^{g,g'}$ and $s_{ij}^{g,g'}$ \eqref{eq:d1}. Then we use equations \eqref{eq:d2} to ensure that $d_{ij}^{g,g'}$ takes the maximum value from $t_{ij}^{g,g'}$ and $s_{ij}^{g,g'}$. Finally, equations \eqref{eq:wd11}--\eqref{eq:wd33} calculate the maximum Wasserstein distances between all sets of combinations of two groups $g$ from $G$ (without repetition, as the order does not matter).

\section{Solution approach} \label{solution}

The proposed Distributionally Fair  P2P optimization model is a bilinear problem. When determining the Wasserstein distance between the different energy groups, the optimization model must calculate the product of the transportation plan (\(\pi_{ij}^{g, g'}\)) and the distance between the transactions (\(d_{ij}^{g, g'}\)). Since both terms are variables, their product introduces a bilinear term to our optimization problem. Therefore, we can establish that constraint \eqref{eq:wd11} is a bilinear constraint. This bilinear constraint introduces non-linearity and non-convexity in our feasible set, which makes it difficult to solve our optimization problem using standard linear algorithms. While the problem size increases, the search for optimality becomes exponentially complex. 

To overcome this issue, we use an alternating algorithm (see Algorithm~\ref{alg:alternating}). The algorithm fixes variables \(d_{ij}^{g, g'}\) transforming them into parameters. Since these variables are now fixed, the resultant model can be directly solved with an LP solver, such as the simplex or interior-point algorithms. The result gives us the current optimal value of the second group of variables \(\pi_{ij}^{g, g'}\). As the next step, the algorithm fixes these second variables \(\pi_{ij}^{g, g'}\) as parameters and runs the optimization model. In this process, it finds the current optimum of the first group of variables \(d_{ij}^{g, g'}\). The algorithm iterates this process until both versions give the same result (i.e., until they converge). 

\begin{algorithm}
\caption{Alternating algorithm}\label{alg:alternating}
\scriptsize
\begin{algorithmic}[1]
\State \textbf{Input:} $\overline{G}$, $c_i$, $e_i$, $\lambda_i^s$, $\lambda_i^b$, $\lambda_i^{ub}$, $\lambda_i^{us}$, $M_i$, $Y_{ij},$$\ \forall i, j \in \mathcal{P}$ 
\State Solve reference model \eqref{eq:ref_obj}--\eqref{eq:ref_non-neg} and get $X$, $U_{i}^{sell}$, $U_{i}^{buy}$, and $e^{c}_{i},  \forall i \in \mathcal{P}$
\State Solve \eqref{eq:revenue} to obtain original revenues ($\gamma_i, \forall i \in \mathcal{P}$)
\State Set deviation level $\varepsilon$, $tol$, and $\overline{iter}$
\State Define initial values $D_{max}^1 \leftarrow \infty$,  $D_{max}^2\leftarrow 0$,  $iter\leftarrow 0$
\While{$|D_{\max}^1 - D_{\max}^2| > tol \land iter \leq \overline{iter}$}
\State Solve  \eqref{eq:transaction} and compute collective transactions ($T_i^g,\ \forall i \in g, \
\forall g \in G $)
\State Solve  \eqref{eq:distance} and compute distances $d^{g, g'}_{ij},\ \forall i \in g, j \in g', \ \forall g, g' \in G: g \neq g'$ (1st variables)
\State Solve \eqref{eq:was1}--\eqref{eq:was2} and compute  $\pi_{ij}^{g, g'},\ \forall i \in g, j \in g', \ \forall g, g' \in G: g \neq g'$ (2nd variables)
\State Solve \eqref{eq:dmax} and identify maximum distance ($D_{\max}^1$)
\State Solve \eqref{eq:fair_obj}--\eqref{eq:wd33} with $\pi_{ij}^{g, g'}$ as fixed parameters and obtain new $X$, and $D_{\max}^2$
\State $iter \leftarrow iter + 1$
\EndWhile
\State $D_{\max}$ $=$ $(D_{\max}^1+D_{\max}^2)/2$
\State \textbf{Output:} $X$, $D_{\max}$
\end{algorithmic}
\end{algorithm}

\section{Case study} \label{case_study}

\begin{figure}[h!]
\centerline{\includegraphics[scale=0.40]{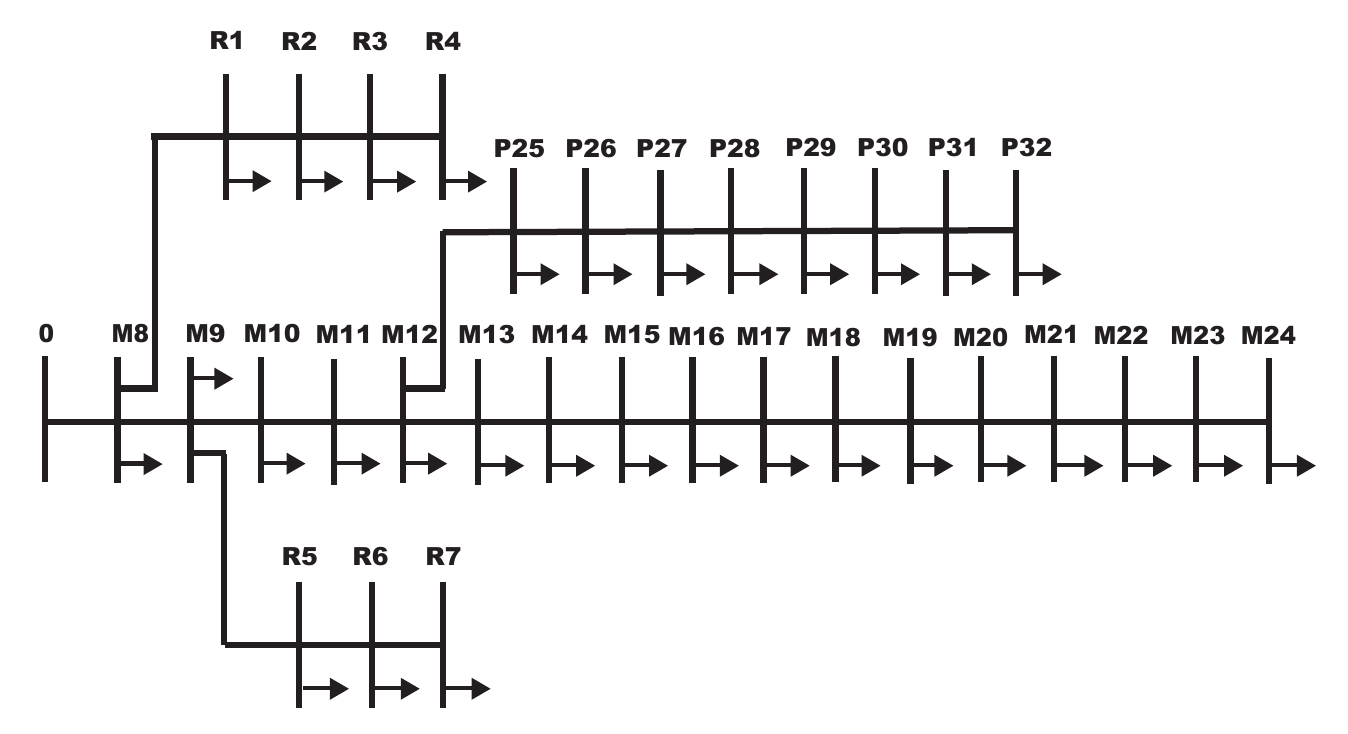}}
\caption{Community's distribution system.}
\label{IEEE_33}
\end{figure}

To evaluate the proposed distributionally fair mechanism for P2P electricity trading, we have used the IEEE 33-bus test system (Fig.~\ref{IEEE_33}), where we enforced the voltage limits $\underline{v} = 0.95$ and $\overline{v} = 1.05$ (in per unit values). We can distinguish households of different energy classes within the distribution grid. Class $\mathscr{R}$ represents households that are energy-rich, class $\mathscr{M}$ represents moderate energy levels, and class $\mathscr{P}$ comprises those households suffering from energy poverty. Within each bus, we connect 50 households, representing an energy community of 1600 households.

\begin{figure}[h!]
\centerline{\includegraphics[scale=0.30]{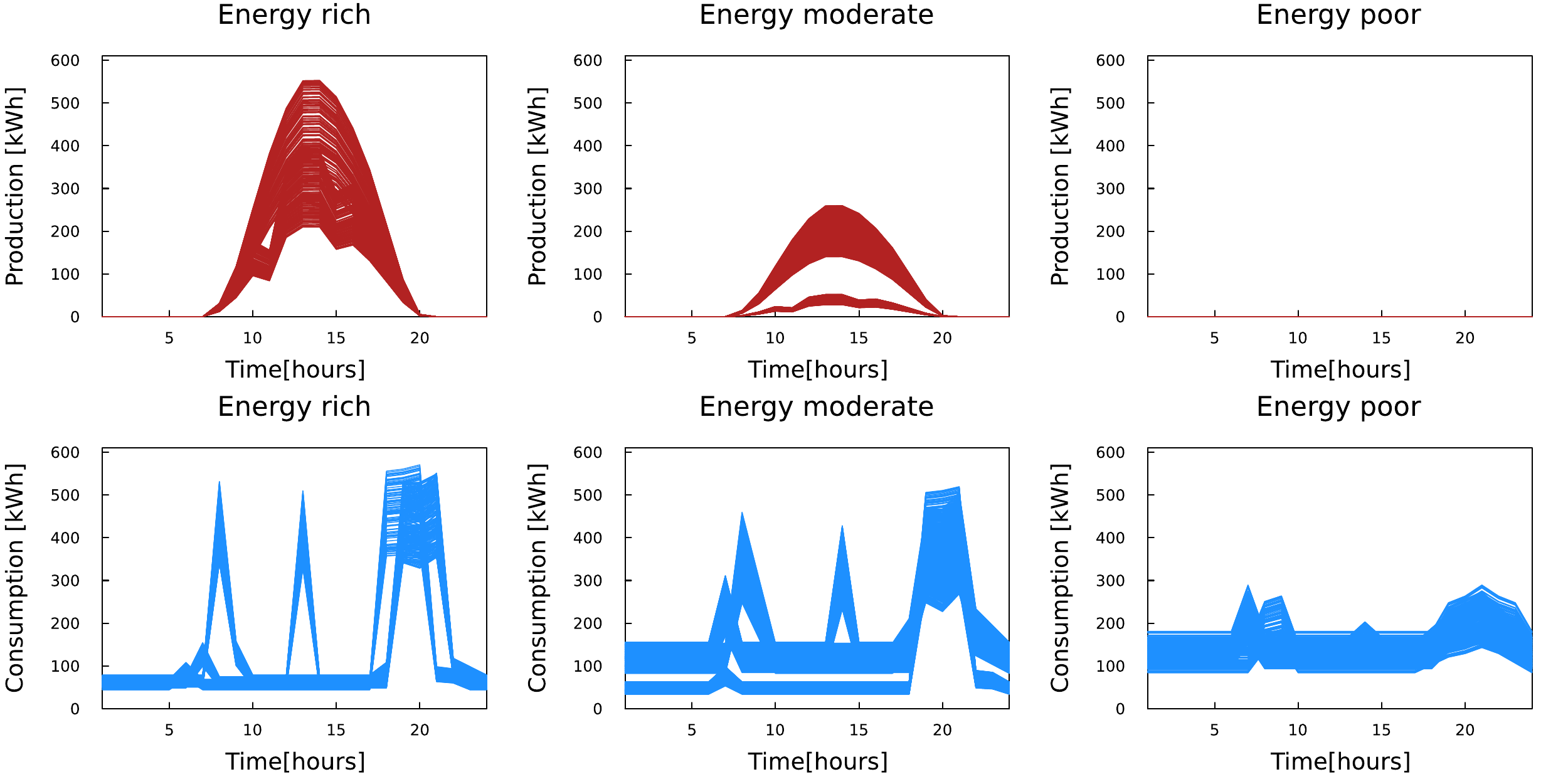}}
\caption{Households' production and consumption profiles.}
\label{scenarios}
\end{figure}

Energy consumption and production profiles are obtained from \citep{IEA_profile}. These profiles have been first normalized to obtain typical consumption and production patterns. Additionally, 10\% of normal noise was added to randomly alter the profiles and make them unique. Then, each profile was scaled by the power peak chosen for each energy group category. Energy-rich households usually have higher power peaks \citep{poverty_peak} due to household size, access to household appliances, and lower payment limitations \citep{poverty_measuring}. Moreover, energy-poor households do not own PV installations generally due to the high investment costs \citep{poverty_PV}. Following these characteristics, we define the following power peaks: 5.1 kW ($\mathscr{R}$), 3.9 kW ($\mathscr{M}$), and 2.1 kW ($\mathscr{P}$). Furthermore, we have considered the following PV ownership shares: 80\% ($\mathscr{R}$), 20\% ($\mathscr{M}$), and 0\%  ($\mathscr{P}$). The resulting consumption and production profiles are represented in Fig.~\ref{scenarios}.
\subsection{Electricity tariffs and bids}

\begin{figure}[h!]
\centerline{\includegraphics[scale=0.28]{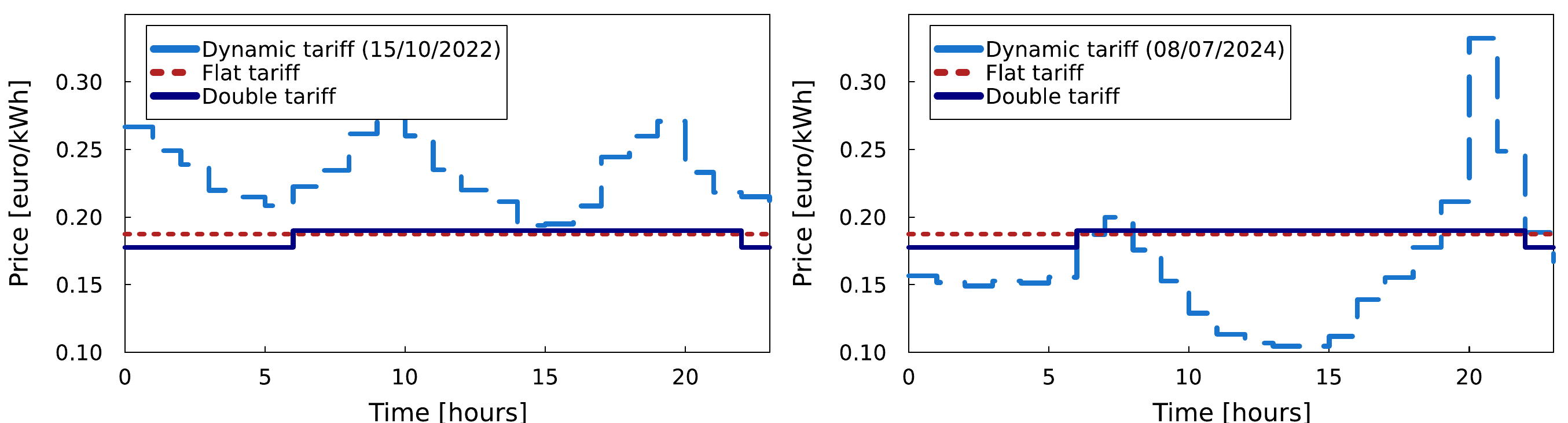}}
\caption{Electricity tariffs on 15/10/2022 (left) and 08/07/2024 (right).}
\label{tariffs}
\end{figure}

Based on the largest Luxembourgish electricity supplier (Enovos), we assume that the utility offers three different electricity tariffs to their customers: flat, double, and dynamic. Within the flat tariff, households pay a price of 0,18736 euro/KWh always \citep{enovos_sell} regardless of the market variations. Moreover, within the double tariff option, households pay a fee of 0,18996 euro/KWh during the day (from 6 am until 10 pm) and 0,17766 euro/KWh during the night (from 10 pm until 6 am) \citep{enovos_sell}. Finally, within the dynamic tariff option, the household is entirely subject to the EPEX clearing price plus a fee of 0,075 euro/KWh for the utility's services  \citep{entsoe_DA}. In Fig.~\ref{tariffs} we represent the different tariffs for the days 15/10/2022 and 08/07/2024. We specifically chose these days because on 15/10/2022 the dynamic price was above the flat tariff and on 08/07/2024 it was below. We define the following tariff distributions within the energy community: 80\% dynamic, 10\% double, and 10\% flat ($\mathscr{R}$); 50\% dynamic, 10\% double, and 40\% flat ($\mathscr{M}$); and 20\% dynamic, 10\% double, and 70\% flat ($\mathscr{P}$). Finally, we consider that the utility buys the surplus energy coming from the PVs at a flat price of 0,1417 euro/KWh \citep{enovos_buy}.

\section{Results} \label{results}

\begin{figure}[h!]
\centerline{\includegraphics[scale=0.25]{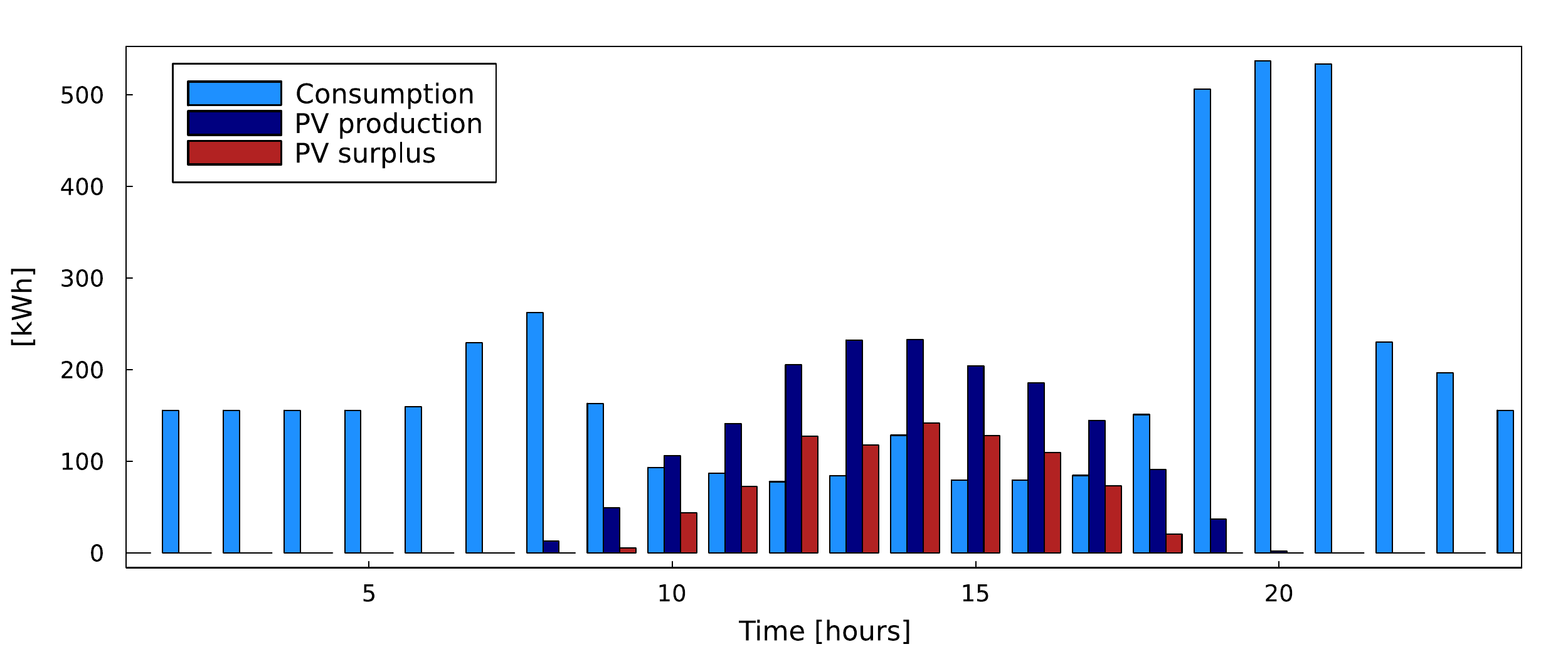}}
\caption{Cummulative PV production, consumption, and surplus.}
\label{surplus}
\end{figure}

The P2P market is not always active. Naturally, when none of the community members has an energy surplus, there will be no transactions within the P2P market. In our case study, the community starts producing PV energy at 8 am (see Fig. \ref{surplus}). However, all the produced energy is self-consumed. It is not until 9 am that the community members start having PV surplus to sell. Then, at 7 pm the PV production is again too small to result in a surplus. Therefore, the P2P market is only active from 9 am until 6 pm. 
\subsection{Base scenario}

In the base scenario, each seller tries to make the highest profit possible from their surplus energy. Therefore, they will always sell their surplus to the highest bidder. As a result, when the energy is traded without further supervision, it is not distributed equally among the different groups of the energy community. The results presented in Table \ref{reference_WD_results} illustrate the evolution of the unfairness level (measured by the maximum Wasserstein distance between the distributions of traded electricity) across different time periods, with a comparison between two dates (15/10/2022 and 08/07/2024). 

\renewcommand{\arraystretch}{0.75}
\begin{table}[h!]
\caption{Evolution of the unfairness level in kWh.}
\centering
\scriptsize
\begin{tabular}{c|ccc|ccc}
\hline
\multicolumn{1}{c|}{\multirow{2}{*}{\textbf{Time}}} & \multicolumn{3}{c|}{\textbf{15/10/2022}} & \multicolumn{3}{c}{\textbf{08/07/2024}} \\ \cline{2-7} 
\multicolumn{1}{c|}{} & \textbf{d\(_\mathscr{RM}\)} & \textbf{d\(_\mathscr{RP}\)} & \textbf{d\(_\mathscr{MP}\)} & \textbf{d\(_\mathscr{RM}\)} & \textbf{d\(_\mathscr{RP}\)} & \textbf{d\(_\mathscr{MP}\)} \\ \hline
{\textbf{00-08:00}} & 0,00 & 0,00 & 0,00 & 0,00 & 0,00 & 0,00 \\ 
\textbf{09:00} & 21,88	& \textbf{22,05	}  & 0,17   & 21,85	& \textbf{22,03}	  &    0,18 \\
\textbf{10:00} & 57,3	& \textbf{95,72}	  & 38,42  & 57,86	& \textbf{95,77}	  &   37,90 \\
\textbf{11:00} & 81,74	& \textbf{102,77 } & 23,94  & 82,78	& \textbf{100,63}	& 23,39 \\
\textbf{12:00} & 141,42	& \textbf{143,72}  & 26,51  & \textbf{145,33}	& 134,67	& 26,73 \\
\textbf{13:00} & 144,36	& \textbf{185,31}  & 41,87  & \textbf{127,02}	& 117,42	& 44,17 \\
\textbf{14:00} & 145,82	& \textbf{213,98}  & 68,16  & \textbf{166,38}	& 136,97	& 67,38 \\
\textbf{15:00} & \textbf{135,09}	& 133,44  & 24,96  & \textbf{137,21}	& 103,57	& 42,41 \\
\textbf{16:00} & 127,46	& \textbf{132,03}  & 23,57  & \textbf{135,65}	& 83,05	  &   55,06 \\
\textbf{17:00} & 97,74	& \textbf{108,29}  & 18,10  &  96,98	& \textbf{109,97}	& 20,99 \\
\textbf{18:00} & 75,17	& \textbf{77,28	}  & 2,12   & 75,1	& \textbf{77,23}	  &   2,14 \\
\textbf{19-23:00} & 0,00 & 0,00    & 0,00   & 0,00 & 0,00 & 0,00 \\ \hline
\textbf{Total} & 1027,98 & \textbf{1214,59} & 267,82 & \textbf{1046,16} & 981,31 & 320,35 \\ \hline
\multicolumn{7}{l}{\textit{*In bold we highlight the maximum Wasserstein distance.}} \\
\hline
\end{tabular}
\label{reference_WD_results}
\end{table}
\renewcommand{\arraystretch}{1.5}

We observe that the unfairness level fluctuates throughout the day, with higher values during periods of increased energy surplus, particularly in the midday hours when solar power generation is at its peak. This is consistent with the effect of PV surpluses (more surplus leads to greater potential for an unfair distribution). Additionally, the utility's tariff also has a clear effect on the fairness level. When dynamic tariffs are higher than flat, peers with dynamic tariffs submit higher bids, thus securing more trades. This is seen on 15/10/2022, when the dynamic tariffs were higher, resulting in greater unfairness between $\mathscr{R}$ and $\mathscr{P}$, since the poor-energy group is mostly subscribed to the flat tariff. In contrast, on 08/07/2024, the dynamic tariffs were lower than the flat tariffs, causing peers with flat tariffs (mainly energy-poor) to submit higher bids and obtain more trades. This shift in tariff dynamics led to a change in the maximum unfairness distance, with the greatest disparity on 08/07/2024 observed between the $\mathscr{R}$ and $\mathscr{M}$, rather than between the $\mathscr{R}$ and $\mathscr{P}$. Thus, the results highlight the significant role that both energy surplus and tariff structures play in the fairness level of the P2P market.

\subsection{A Distributionally Fair P2P market}

\begin{figure}[h!]
\centerline{\includegraphics[scale=0.35]{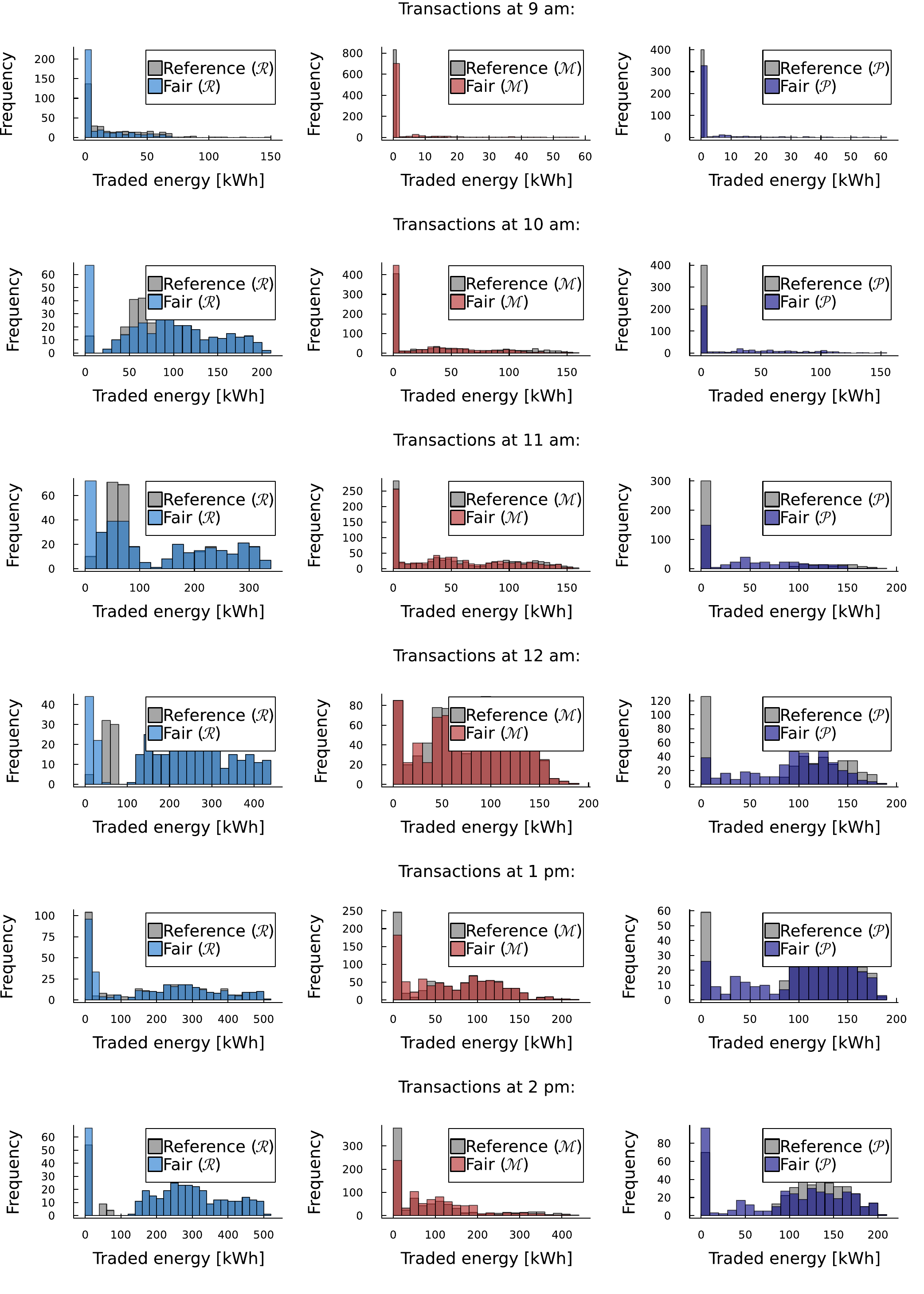}}
\caption{Trade distributions in reference and fair ($\varepsilon=100\%$) scenarios (08/07/2024 - I).}
\label{histogram1}
\end{figure}

\begin{figure}[h!]
\centerline{\includegraphics[scale=0.35]{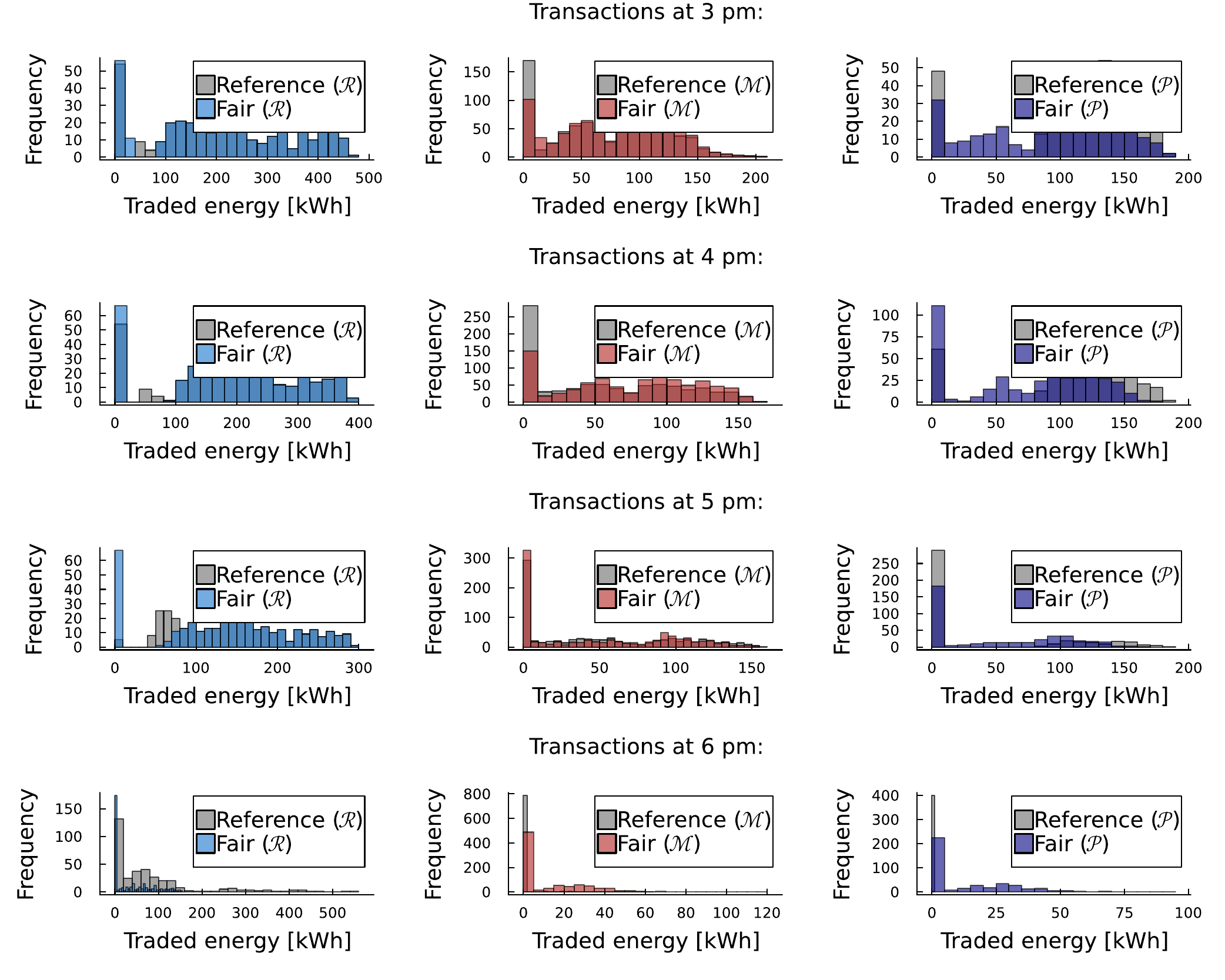}}
\caption{Trade distributions in reference and fair ($\varepsilon=100\%$) scenarios (08/07/2024 - II).}
\label{histogram2}
\end{figure}

When applying the Distributionally Fair mechanism for P2P trading, the goal is not to increase the revenue of the sellers but to decrease the unfairness level between the different groups. The model redistributes the trades among the peers to minimize the maximum Wasserstein distance between the trade distributions associated with each group (see Fig. \ref{histogram1} and \ref{histogram2}). However, the model can apply fairness with a limitation, given by a user-defined sacrifice level \(\varepsilon\), which corresponds to the maximum percentage of revenue sacrifice with respect to the reference scenario that the market is willing to tolerate.

Table \ref{fair_results} identifies the maximum unfairness level for different sacrifice values. We have used the alternating algorithm for this identification as explained in Section \ref{solution}. For the lowest sacrifice level ($\epsilon = 1\%$), the initialization of the algorithm was based on the transaction matrix $X$ that is obtained from the reference P2P model. Then, for higher sacrifice levels, the initialization relied on the transaction matrix obtained from the previously solved distributionally fair P2P model, with the closest preceding sacrifice level. For instance, the results for $\epsilon = 2\%$ were derived using the transaction matrix from $\epsilon = 1\%$, the results of $\epsilon = 5\%$ were obtained by initializing the algorithm with the energy transactions corresponding to $\epsilon = 2\%$, and so on. Additionally, the tolerance level \emph{tol} was set to 0.01 and the maximum iterations $\overline{iter}$ to 15.

\renewcommand{\arraystretch}{0.75}
\begin{table}[h!]
\caption{Unfairness level (kWh) using the Distributionally Fair P2P market (08/07/2024).}
\resizebox{\textwidth}{!}{
\centering
\begin{tabular}{c|c|cccccccc}
\hline
\multicolumn{1}{c|}{\multirow{2}{*}{\textbf{Time}}} & \multicolumn{1}{c|}{\multirow{2}{*}{\textbf{Ref.}}} & \multicolumn{8}{c}{\textbf{Distributionally Fair model ($\epsilon$)}} \\ \cline{3-10} 
\multicolumn{1}{c|}{} & \multicolumn{1}{c|}{} & \textbf{1\%} & \textbf{2\%} & \textbf{5\%} & \textbf{10\%}  & \textbf{20\%} & \textbf{50\%} & \textbf{70\%}& \textbf{100\%}  \\ \hline
\textbf{09:00}    & 22,03	& 21,62  &  21,06  &  19,70    & 17,10   & 12,42 & \textbf{7,57}   & 7,57  & 7,57  \\
\textbf{10:00}    & 95,77	& 93,15  &  87,98  &  75,32    & 64,42   & \textbf{58,03} & 58,03  & 58,03 & 58,03 \\
\textbf{11:00}    & 100,63	& 94,95  &  88,22  &  \textbf{74,80}    & 74,80   & 74,80 & 74,80  & 74,80 & 74,80 \\
\textbf{12:00}    & 145,33	& 141,08 &  135,24 &  \textbf{135,23}   & 135,23  & 135,23& 135,23 & 135,23 & 135,23 \\
\textbf{13:00}    & 127,02	& 121,80 &  \textbf{117,84} &  117,84   & 117,84  & 117,84& 117,84 & 117,84 & 117,84 \\
\textbf{14:00}    & 166,38	& 162,77 &  161,70 &  159,10   & 154,14  & \textbf{153,96} & 153,96 & 153,96 & 153,96 \\
\textbf{15:00}    & 137,21	& 133,51 &  133,49 &  133,49   & 133,49  & 133,49& 133,49 & \textbf{132,13} & 132,13 \\
\textbf{16:00}    & 135,65	& 132,54 &  131,02 &  128,51   & 123,74  & \textbf{117,06} & 117,06 & 117,06 & 117,06\\
\textbf{17:00}    & 109,97	& 105,35 &  98,71  &  89,41    & \textbf{87,44}   & 87,44 & 87,44  & 87,44 & 87,44 \\
\textbf{18:00}    & 77,23	& 75,74  &  72,87  &  68,04    & 58,93   & 42,53 & \textbf{23,10}  & 23,10 & 23,10 \\ \hline
\textbf{tot}      & 1.117,22& 1.082,51 &1048,15&1.001,44 &	967.13   & 932,80& 908,52 & \textbf{907,16} & 907,16\\
\hline
\multicolumn{10}{l}{\textit{*In bold we highlight the fairness limit and its corresponding sacrifice level. }} \\
\hline
\end{tabular}
}
\label{fair_results}
\end{table}
\renewcommand{\arraystretch}{1.5}

We learn that, as expected, the higher the sacrifice level $\epsilon$ the market tolerates, the lower the level of unfairness between the different groups in the community. The distributionally fair P2P model decreases the unfairness level up to 70.1\% (at 18:00) and by 23\% on average. Nonetheless, there is a limit to the achievable fairness in the problem. At one point, it does not matter if we further increase the sacrifice level, the unfairness level will not decrease any more. For some hours, we see that this limit is reached very quickly (at 13:00 with $\varepsilon=2\%$). Conversely, in other hours fairness keeps improving while increasing the sacrifice level (at 15:00 until $\varepsilon=70\%$). When a small number of energy-rich households generate surplus energy, which must be purchased by the larger number of moderate and energy-poor households, achieving a perfectly equitable trading scenario becomes infeasible. Since surplus energy is typically less than the overall deficit, energy-rich households are assured of finding buyers. However, energy-moderate and poor households must compete for these trades, leading to a situation where energy-rich households dominate trade participation. Consequently, even when we allow for a full sacrifice level, the unfairness can still go up to 153.96 kWh (see Table \ref{fair_results}).

\subsection{Addition of a non-profit community PV plant}

In an attempt to further decrease the unfairness level within the energy community, a nonprofit shareable PV plant can be incorporated  (connected to node 12 in our case study). The production pattern of the PV plant was obtained from \citep{IEA_profile}. We normalized it and then scaled it for different installed capacities (5 kW, 10 kW, 15 kW, and 20 kW). The generation of the community PV plant is offered within the P2P market, aiming to decrease the unfairness level as much as possible. Table~\ref{fair_PV} illustrates that as the installed capacity of the PV plant increases, the unfairness level declines. The introduction of a 20 kW PV plant can reduce the unfairness level by 51.95\% on average, with a nearly complete elimination achieved at certain hours. 

\renewcommand{\arraystretch}{0.75}
\begin{table}[h!]
\caption{Unfairness level (kWh) with a non-profit community PV plant (08/07/2024).}
\centering
\scriptsize
\begin{tabular}{c|c|ccccc}
\hline
\multicolumn{1}{c|}{\multirow{2}{*}{\textbf{Time}}} & \multicolumn{1}{c|}{\multirow{2}{*}{\textbf{Ref.}}} & \multicolumn{5}{c}{\textbf{Distributionally Fair model (\(\varepsilon\)=1.0) }} \\ \cline{3-7} 
\multicolumn{1}{c|}{} & \multicolumn{1}{c|}{} & \textbf{0 kW} & \textbf{5 kW} & \textbf{10 kW} & \textbf{15 kW}  & \textbf{20 kW}\\ \hline
\textbf{09:00}    & 22,03	& 7,57   &  4,54      & 1,50   & 0,02   & 0,01 \\
\textbf{10:00}    & 95,77	& 58,03  &  51,51     & 44,98  & 38,46  & 31,94 \\
\textbf{11:00}    & 100,63	& 74,80  &  69,08     & 63,37  & 57,65  & 51,94 \\
\textbf{12:00}    & 145,33	& 135,30 &  122,59    & 109,94 & 97,46  & 84,67 \\
\textbf{13:00}    & 127,02	& 118,99 &  104,43    & 89,85  & 75,47  & 61,48 \\
\textbf{14:00}    & 166,38	& 153,96 &  139,67    & 125,39 & 111,10 & 96,82 \\
\textbf{15:00}    & 137,21	& 130,13 &  120,52    & 108,99 & 98,96  & 87,42 \\
\textbf{16:00}    & 135,65	& 115,59 &  104,40    & 93,71  & 81,81  & 69,93 \\
\textbf{17:00}    & 109,97	& 87,44  &  78,56     & 69,67  & 60,78  & 51,90 \\
\textbf{18:00}    & 77,23	& 23,10  &  17,51     & 11,92  & 6,33   & 0,75 \\ \hline
\textbf{tot}      & \textbf{1.117,22} & \textbf{904,90}& \textbf{812,80} &  \textbf{719,32} & \textbf{628,04} & \textbf{536,84}\\
\hline
\end{tabular}
\label{fair_PV}
\end{table}
\renewcommand{\arraystretch}{1.5}

\subsection{Computational performance analysis}

The reference P2P model and the distributionally fair P2P model differ fundamentally in their approach and computational complexity. The reference model facilitates standard P2P electricity trading, relying on a linear optimization framework. On the other hand, the distributionally fair P2P model integrates equity principles into the trading process. This inclusion transforms the problem into a challenging bilinear optimization problem, requiring an alternating algorithm to solve it, which involves multiple iterations and a higher computational effort. 

The results in Table~\ref{computation} highlight the difference in computational requirements between the two models. The reference model, solved using an LP algorithm via Gurobi, maintains an average computational cost of 2 minutes and 5 seconds, demonstrating its efficiency with minimal computation time for each time slot. Conversely, the distributionally fair P2P model exhibits varying computational demands depending on the time slot and sacrifice level, ranging from 1 minute and 25 seconds to 25 minutes and 54 seconds. This variation arises from the increased complexity of ensuring a fairer electricity trade exchange under conditions of higher energy demands and generation levels. Additionally, the number of iterations required fluctuates with the hour, spanning from 2 to 12 iterations. These findings highlight the extra layer of complexity that arises from constraining the outcomes of a P2P electricity market with some measure of \emph{distributional fairness}. 

\renewcommand{\arraystretch}{0.75}
\begin{table}[h!]
\scriptsize
\caption{Measured computational times (in seconds).}
\centering
\begin{tabular}{c|c|cccccccccc}
\hline
\multicolumn{1}{c|}{\multirow{2}{*}{\textbf{Time}}} & \multicolumn{1}{c|}{\multirow{2}{*}{\textbf{\begin{tabular}[c]{@{}c@{}}Ref. \\ \end{tabular}}}} & \multicolumn{8}{c}{\textbf{Distributionally Fair model ($\varepsilon$)}} \\ \cline{3-10} 
\multicolumn{1}{c|}{} & \multicolumn{1}{c|}{} & \multicolumn{1}{c}{\textbf{1\%}} & \multicolumn{1}{c}{\textbf{2\%}} & \textbf{5\%} & \multicolumn{1}{c}{\textbf{10\%}} & \multicolumn{1}{c}{\textbf{20\%}} &\multicolumn{1}{c}{\textbf{50\%}} &\multicolumn{1}{c}{\textbf{70\%}} &\multicolumn{1}{c}{\textbf{100\%}} \\ \hline
\textbf{09:00} & 125 & 232 & 226 & 248 & 238 & 256 & 420 & 126 & 98 \\
\textbf{10:00} & 123 & 235 & 213 & 247 & 221 & 244 & 105 & 118 & 96 \\
\textbf{11:00} & 134 & 244 & 684 & 521 & 120 & 99 & 127 & 101 & 98 \\
\textbf{12:00} & 115 & 247 & 1492 & 115 & 107 & 126 & 103 & 108 & 127 \\
\textbf{13:00} & 122 & 1191 & 935 & 111 & 107 & 129 & 125 & 108 & 126 \\
\textbf{14:00} & 142 & 242 & 272 & 252 & 241 & 252 & 135 & 129 & 147 \\
\textbf{15:00} & 116 & 472 & 265 & 165 & 145 & 147 & 124 & 742 & 151 \\
\textbf{16:00} & 127 & 387 & 254 & 276 & 256 & 1554 & 129 & 112 & 122 \\
\textbf{17:00} & 132 & 224 & 281 & 396 & 423 & 156 & 151 & 121 & 110 \\
\textbf{18:00} & 116 & 186 & 190 & 184 & 180 & 198 & 411 & 85 & 86 \\ \hline
\textbf{Average} & \textbf{125} & \textbf{366} & \textbf{481} & \textbf{252} & \textbf{204} & \textbf{316} & \textbf{183} & \textbf{175} & \textbf{116} \\\hline
\multicolumn{9}{l}{\textit{*Ref: Reference P2P model}} \\ \hline
\end{tabular}
\label{computation}
\end{table}
\renewcommand{\arraystretch}{1.5}

\section{Conclusions} \label{conclusions}

Within this article, we considered a semi-decentralized P2P market architecture with a community manager using the IEEE 33-bus case study and 1600 peers. We separate the peers into three groups according to their energy poverty level: energy-rich, energy-moderate, and energy-poor. We demonstrate that without further supervision, the electricity trades are not distributed equally among the different groups since peers behave selfishly. 

To bring equity into the P2P electricity trades, we have proposed and analyzed a distributionally fair P2P optimization model. We solved the original NP-hard problem using an alternating algorithm and demonstrated that the proposed technique can better distribute the electricity trades to reduce the unfairness level. The proposed model can identify the discriminated group and take appropriate action. Unlike other models, it can identify a different group depending on the specific conditions of the moment. Therefore, it will help all groups when needed. Nonetheless, complete fairness cannot be achieved even with a full sacrifice level. This is primarily because the energy-rich community generates most of the surplus, allowing them to secure more trades than other groups. We demonstrate that we can further decrease the unfairness level by including a nonprofit community PV plant.

For future work, there are numerous opportunities to refine and expand upon the proposed model. One key limitation of the current solution algorithm is its susceptibility to converging to local minima. Although we have proven that the distributionally fair P2P model can increase fairness, we cannot guarantee that we identify the global optimum. Bilinear problems can be reformulated as mixed-integer linear programming optimization models. The alternating algorithm's task of fixing one variable at a time can be accomplished by introducing an auxiliary binary variable. This approach can guarantee global optimum. However, the size of this auxiliary binary variable depends on the number of peers involved, which leads to a significantly more complex final optimization model in practice. Therefore, further research is needed to address this challenge and enhance the formulation for more efficient solving in real-world applications.

Finally, while this study defined groups based on energy poverty levels, the model can be adapted to define these groups using other criteria, such as the type of tariff households have, the amount of PV ownership, or even their consumption patterns. This flexibility opens up further applications of the model in different contexts, allowing for its use in various markets and regions where equity is a concern.

\section*{Acknowledgments}
E. R. Irusta gratefully acknowledges the financial support of Creos Luxembourg S.A. under the research project FlexBeAn; the Luxembourg National Research Fund (FNR) with grant reference 17742284, and PayPal, PEARL grant reference 13342933/GF. For the purpose of open access, and in fulfillment of the obligations arising from the grant agreement, the author has applied a Creative Commons Attribution 4.0 International (CC BY 4.0) license to any Author Accepted Manuscript version arising from this submission. The work of J. M. Morales was supported by the Spanish Ministry of Science and Innovation (AEI/10.13039/501100011033) through project PID2023-148291NB-I00.

\bibliographystyle{elsarticle-harv} 
\bibliography{ref}

\end{document}